\documentclass[12pt]{article}
\usepackage{amsmath}
\usepackage{amssymb}
\usepackage{dsfont}
\usepackage{cite}
\newsymbol \blackbox 1004
\newcommand{\eh}{\hfill}\newlength{\sperr}

\newenvironment{proof}{{\settowidth{\sperr}{\bf\rm
Proof}%
\par\addvspace{0.3cm}\noindent\parbox[t]{1.3\sperr}
{\bf\rm P\eh r\eh o\eh o\eh f\eh }%
}}{\nopagebreak\mbox{}
$\blackbox$\par\addvspace{0.3cm}}

\def\nn{\nonumber}

\def\a{\alpha}

\def\Lam{\Lambda}
\def\s{\sigma}
\def\la{\lambda}

\def\wt{\widetilde}
\def\ov{\overline}

\def\BC{{\mathbb C}}

\def\BN{{\mathbb N}}

\def\cli{{\mathcal I}}

\def\cls{\mathcal{S}}

\def\diag{\mathrm{diag}}
\newcommand{\E}{\mathrm{e}}
\newcommand{\I}{\mathrm{i}}

\newtheorem{Pa}{Paper}[section]
\newtheorem{Tm}[Pa]{{\bf Theorem}}

\newtheorem{Rk}[Pa]{{\bf Remark}}

\newtheorem{Dn}[Pa]{{\bf Definition}}

\newtheorem{Pn}[Pa]{{\bf Proposition}}

\newenvironment{dedication}
        {\vspace{1ex}\begin{quotation}\begin{center}\begin{em}}    
        {\par\end{em}\end{center}\end{quotation}}

\title{GBDT of discrete skew-selfadjoint Dirac systems and explicit solutions
of the corresponding non-stationary problems}

\author{Alexander Sakhnovich}

\date{}
\begin{document}
\maketitle

\begin{dedication}
Dedicated  to Rien Kaashoek    on the occasion of his 80th anniversary
\end{dedication}

\begin{abstract}    Generalized B\"acklund-Darboux transformations (GBDTs) of
discrete skew-selfadjoint Dirac systems have been successfully used for
explicit solving of direct and inverse problems of Weyl-Titchmarsh theory. During explicit solving
 of direct and inverse problems, we considered GBDT of the trivial initial
 systems. However, GBDT of arbitrary discrete skew-selfadjoint Dirac systems
 is important as well and we introduce these transformations in the present paper.
 The obtained results are applied to the construction of explicit solutions
 of the interesting related  non-stationary systems.
\end{abstract}

{MSC(2010):  34A05, 39A06, 39A12.}

\vspace{0.2em}

Keywords:  {\it Discrete skew-selfadjoint Dirac system, generalized B\"acklund-Darboux transformation, fundamental solution, non-stationary system, explicit solution.}

\section{Introduction}\label{Intro}
\setcounter{equation}{0}
We consider a discrete skew-selfadjoint Dirac system
\begin{equation} \label{R1}
y_{k+1}(z)=\left(I_m+ \frac{\I}{ z}
C_k\right)
y_k(z),  \quad C_k=U_k^{*}jU_k \quad \left( k \in \cli \right),
\end{equation}
where $I_m$  is the $m\times m$ identity matrix, $U_k$ are $m \times m$ unitary matrices,
\begin{align} &   \label{I2}
j = \left[
\begin{array}{cc}
I_{m_1} & 0 \\ 0 & -I_{m_2}
\end{array}
\right]  \quad (m_1, m_2 \in \BN,  \quad  m_1+m_2=m),
 \end{align} 
 and $\cli$ is either $\BN_0$ or the set $\{k\in \BN_0: \, 0 \leq k < N<\infty\}$ $(N\in \BN)$. 
 Here, as usual, $\BN$ denotes the set of natural numbers and $\BN_0=0 \cup \BN$.

 The relations
 \begin{align} \label{R2-}&
C_k=C_k^*, \quad C_k^2=I_m
\end{align}
 are immediate from the second equality in \eqref{R1}.

This paper is a certain prolongation of the papers \cite{KS, FKKS}, where  direct and inverse problems
for  Dirac systems \eqref{R1} have been solved explicitly and explicit solutions  of the
isotropic Heisenberg magnet model have been constructed. We would like to mention also
the earlier papers on the cases of the continuous Dirac systems (see, e.g., \cite{GKS2, GKS6}).

The GBDT version of the B\"acklund-Darboux transformations have been used in \cite{GKS2, GKS6, KS, FKKS}.
B\"acklund-Darboux transformations and related commutation methods (see, e.g., \cite{Ci, D, Ge, GeT,  MS, SaSaR, ZM}
and numerous references therein) are well-known tools in the spectral theory and in the construction of explicit
solutions. In particular, the generalized B\"acklund-Darboux transformations 
(i.e., the  GBDT version of the B\"acklund-Darboux transformations) were introduced in \cite{SaA1} and developed further in a series of papers
(see \cite{SaSaR} for details). 

Whereas GBDTs of the trivial initial
systems have been used in \cite{GKS2, GKS6, KS, FKKS} (in particular, initial systems \eqref{R1} where $C_k\equiv j$
have been considered in \cite{KS, FKKS}), the case of an arbitrary initial  discrete Dirac system \eqref{R1} is considered here.
In Section \ref{GBDT1}, we introduce GBDT, construct the so called Darboux matrix and give representation of the fundamental
solution of the transformed system (see Theorem \ref{Fund}). Note that an explicit representation of the fundamental
solutions of the transformed systems in terms of the solutions of the initial systems is one of the main features and advantages
of the Darboux transformations.

One of the recent developments of the GBDT theory is connected with its application to the construction
of explicit solutions of important dynamical systems (see, e.g., \cite{FKRS, SaA17}).
In Section 3 of this article, we use the same approach in order to construct explicit solutions of
the non-stationary systems corresponding to the systems \eqref{R1}.

 In the paper, $\BC$ stands for the complex plane
and $\BC_+$ stands for the open upper halfplane.
The notation $\s(\a)$ stands for the spectrum of the matrix $\a$ and the notation $\diag\{d_1, d_2, \ldots\}$
stands for the block diagonal matrix with the blocks $d_1, d_2, \ldots$ on the main diagonal.

\section{GBDT of discrete skew-selfadjoint \\ Dirac systems}\label{GBDT1}
\setcounter{equation}{0}
Each GBDT of the initial system \eqref{R1} is determined by some
triple $\{\a, S_0, \Lam_0\}$ of the $n \times n$ matrices $\a$ and $S_0=S_0^*$  and the $n \times m$ matrix $\Lam_0$ ($n \in \BN$)
such that 
\begin{align} \label{R2}&
\a S_0-S_0\a^*=\I\Lam_0 \Lam_0^*.
\end{align}
The initial skew-selfadjoint Dirac system has the form \eqref{R1} and the transformed (i.e., GBDT-transformed)
system has the form
\begin{equation} \label{R1'}
\wt y_{k+1}(z)=\left(I_m+ \frac{\I}{ z}
\wt C_k\right)
\wt y_k(z) \quad \left( k \in \cli \right),
\end{equation}
where the {\it potential} $\{\wt C_k\}$ $(k\in \cli)$
is given  by the relations
\begin{align}\label{R3}
&\Lam_{k+1}= \Lam_{k}+\I \a^{-1} \Lam_{k} C_k, 
 \\
&
S_{k+1}=S_k+ \a^{-1} S_k (\a^*)^{-1}+ \a^{-1} \Lam_{k} C_k \Lam_{k}^*
(\a^*)^{-1},\label{R4} \\
&\label{R5}
\wt C_k=C_k+ \Lam_k^* S_k^{-1} \Lam_k - \Lam_{k+1}^* S_{k+1}^{-1} \Lam_{k+1}, \quad k \in \cli.
\end{align}
Here and further in the text we assume that 
\begin{equation} \label{R6}
\det \a \not=0, 
\end{equation}
and suppose additionally in \eqref{R5} that
\begin{equation} \label{R7}
 \det S_k\not=0 
\end{equation}
for $k\in \BN_0$ or for $0\leq k \leq N$ depending on the choice of the interval $\cli$, on which the Dirac system is considered.

Similar to the proof of \cite[(3.7)]{FKKS}, using the equality $C_k^2=I_m$  from \eqref{R2-} and relations \eqref{R2}--\eqref{R4} 
one  easily proves by induction
that
\begin{align} \label{R8}&
\a S_k-S_k \a^*=\I\Lam_k \Lam_k^*.
\end{align}

\begin{Rk} Clearly, $S_k=S_k^*$ and $\wt C_k=\wt C_k^*$. Further in the text, in Proposition \ref{PnCk} we show that $\wt C_k^2=I_m$. 
In Theorem \ref{PnCk+}, we show that under conditions $S_0>0$ and $0, \I \not\in \s(\a)$ $(\s(\a)$ is the spectrum of $\a)$
we have $S_k>0$ and
\begin{align} \label{R1+}&
\wt C_k= \wt U_k^{*}j\wt U_k ,
\end{align}
where
$j$ is given in \eqref{I2} and the matrices $\wt U_k$ are unitary. The equality 
\eqref{R1+} means that the transformed system \eqref{R1'}
is again a skew-selfadjoint Dirac system in the sense of the definition \eqref{R1}.
Before  Theorem \ref{PnCk+} we consider the GBDT-transformed system
\eqref{R1'} without the requirement \eqref{R1+}.
\end{Rk}
The fundamental solutions of \eqref{R1} and \eqref{R1'} are denoted by 
$w(k, z)$ and $\wt w(k,z)$, respectively,
and are normalized  by the conditions
\begin{align} \label{R9}&
\wt w(0, z)= w(0,z)=I_m.
\end{align}
In other words, $y_k=w(k, z)$ and $\wt y_k= \wt w(k, z)$ are $m\times m$ matrix
solutions of the initial and transformed systems, respectively, which
satisfy the initial conditions \eqref{R9}.
The so called Darboux matrix corresponding to the transformation of the system
\eqref{R1} into \eqref{R1'} is given by the transfer matrix function $w_{\a}$
in Lev Sakhnovich form:
\begin{align} \label{R10}&
w_{\a}(k, z)=I_m-\I \Lam_k^*S_k^{-1}(\a- z I_n)^{-1}\Lam_k.
\end{align}
See \cite{SaL1, SaL2} as well as \cite{SaSaR} and further references therein for
the notion and properties of this transfer matrix function. The statement that 
the Darboux matrix has the form \eqref{R10} may be formulated as the following
theorem.
\begin{Tm}\label{Fund} Let the initial Dirac system \eqref{R1} and a
triple $\{\a, S_0, \Lam_0\}$, which satisfies the relations  \eqref{R2}, \eqref{R6}, \eqref{R7} and $S_0=S_0^*$, be given. Then, the fundamental solution $w$ of the
initial system and fundamental solution $\wt w$ of the transformed system
\eqref{R1'} $($determined by the triple $\{\a, S_0, \Lam_0\}$ via relations
\eqref{R3}--\eqref{R5}$)$ satisfy
 the equality
\begin{align} \label{R11}&
\wt w(k, z)= w_{\a}(k,-z)w(k,z)w_{\a}(0,-z)^{-1} \quad (k \geq 0),
\end{align} 
where $w_{\a}$ has the form \eqref{R10}. 
\end{Tm} 
\begin{proof}.
The following equality is crucial for our proof
\begin{equation} \label{R12}
w_{ \alpha }(k+1, z )\left(I_{m}-
\frac{\I}{z}C_k\right)=\left(I_{m}- \frac{\I}{z}\wt C_k\right)w_{ \alpha }(k,
z).
\end{equation}
(It is easy to see that the important formula \cite[(3.16]{FKKS} is a particular case of \eqref{R12}.)
In order to prove \eqref{R12}, note that according to (\ref{R10}) formula (\ref{R12}) is equivalent to the formula
\begin{align}
\frac{1}{z}\left(\wt C_k-C_k\right)=&\left(I_{m}- \frac{\I}{z}\wt C_k\right) \Lam_k^{*}
S_k^{-1} ( z I_{n}- \alpha )^{-1} \Lam_k
\nonumber\\ & -
\Lam_{k+1}^{*} S_{k+1}^{-1}
(z I_{n}- \alpha )^{-1} \Lam_{k+1}\left(I_{m}-
\frac{\I}{z} C_k\right).\label{R13}
\end{align}
Using the Taylor expansion of $(z I_{n}- \alpha )^{-1}$ at
infinity we see that (\ref{R13}) is in  turn equivalent to
the set of equalities:
\begin{align}&
\wt C_k-C_k=\Lam_k^{*} S_k^{-1}
\Lam_k - \Lam_{k+1}^{*} S_{k+1}^{-1}
\Lam_{k+1},\label{R14}\\
& \Lam_{k+1}^{*} S_{k+1}^{-1} \a^p
\Lam_{k+1}-\I \Lam_{k+1}^{*} S_{k+1}^{-1} \a^{p-1}
\Lam_{k+1} C_k\nonumber\\ &
= \Lam_k^{*} S_{k}^{-1} \a^p \Lam_k-\I
\wt C_k \Lam_k^{*} S_{k}^{-1} \a^{p-1} \Lam_k \qquad
(p>0).\label{R15}
\end{align}
Equality (\ref{R14}) is equivalent to (\ref{R5})
and it remains to prove \eqref{R15}. From \eqref{R3}, taking into account $C_k^2=I_m$ we have 
\begin{align}&
\a \Lam_{k+1} -\I
\Lam_{k+1} C_k=\a \Lam_{k+1} -\I
\Lam_k C_k+ \a^{-1} \Lam_k=\a \Lam_k + \a^{-1} \Lam_k. 
\label{R16}
\end{align}
Substituting \eqref{R16} into the left hand side of \eqref{R15} and using simple transformations,
we  rewrite \eqref{R15} in the form 
\begin{align} \nn &
Z_k \a^{p-2} \Lam_k=0, \quad Z_k:=\Lam_{k+1}^{*} S_{k+1}^{-1} (\a^2+I_{n})-
\Lam_k^{*} S_{k}^{-1} \a^2+\I \wt C_k \Lam_k^{*}
S_{k}^{-1} \a.
\end{align}
Therefore, in order to prove \eqref{R15} (and so to prove \eqref{R12}) it suffices  to show that
\begin{align} \label{R17} &
\Lam_{k+1}^{*} S_{k+1}^{-1} (\a^2+I_{n})=
\Lam_k^{*} S_{k}^{-1} \a^2-\I \wt C_k \Lam_k^{*}
S_{k}^{-1} \a,
\end{align}
that is, $Z_k=0$. Relation \eqref{R17} is of interest in itself, since it is
an analogue of \eqref{R3} (more precisely of the relation adjoint to \eqref{R3})
when $\Lam_{r}^{*}S_r^{-1}$ is taken instead of $\Lam_r^*$. Such analogues
are useful in continuous and discrete GBDT as well as in the construction of explicit
solutions of dynamical systems  (see, e.g. \cite{FKRS, SaA17, SaSaR}
and references therein).

Taking into account  \eqref{R5} and \eqref{R3}, we rewrite \eqref{R17} in the form
\begin{align}\nn
&\Lam_{k+1}^{*} S_{k+1}^{-1} (\a^2+I_{n})-
\Lam_k^{*} S_{k}^{-1} \a^2 \\ & \nn
\quad +\I
\big(C_k+ \Lam_k^* S_k^{-1} \Lam_k - \Lam_{k+1}^* S_{k+1}^{-1} \Lam_{k+1}\big)
 \Lam_k^{*}
S_{k}^{-1} \a
\\ & \nn
= \Lam_{k+1}^{*} S_{k+1}^{-1} (\a^2+I_{n})- \I \Lam_{k+1}^* S_{k+1}^{-1} (\Lam_k+\I \a^{-1} \Lam_k C_k)\Lam_k^{*}
S_{k}^{-1} \a
\\ & \quad -
\Lam_k^{*} S_{k}^{-1} \a^2 +\I
\big(C_k+ \Lam_k^* S_k^{-1} \Lam_k \big)
 \Lam_k^{*}
S_{k}^{-1} \a=0.
\label{R18}
\end{align}
Since \eqref{R8} yields
\begin{align}
\I \Lam_k \Lam_k^{*}
S_{k}^{-1}= \a - S_k \a^* S_k^{-1},
\label{R19}
\end{align}
we rewrite the third line in \eqref{R18} and see that \eqref{R18} (i.e., also \eqref{R17}) is equivalent to
\begin{align} & \nn
\Lam_{k+1}^{*} S_{k+1}^{-1} (I_{n}+ \a^{-1} \Lam_k C_k\Lam_k^{*}
S_{k}^{-1} \a + S_k \a^* S_k^{-1}\a)
\\ &  -
\Lam_k^{*} S_{k}^{-1} \a^2 +\I
\big(C_k+ \Lam_k^* S_k^{-1} \Lam_k \big)
 \Lam_k^{*}
S_{k}^{-1} \a=0.
\label{R20}
\end{align}
Formula \eqref{R4} implies that
$$(I_{n}+ \a^{-1} \Lam_k C_k\Lam_k^{*}
S_{k}^{-1} \a + S_k \a^* S_k^{-1}\a)=S_{k+1}\a^* S_k^{-1}\a .$$
Hence, \eqref{R20} is equivalent to
\begin{align} & 
\Lam_{k+1}^{*} \a^* S_k^{-1}\a  -
\Lam_k^{*} S_{k}^{-1} \a^2 +\I
\big(C_k+ \Lam_k^* S_k^{-1} \Lam_k \big)
 \Lam_k^{*}
S_{k}^{-1} \a=0.
\label{R21}
\end{align}
Using again \eqref{R19}, we rewrite \eqref{R21} as
\begin{align} & 
\big(\Lam_{k+1}^{*}-\Lam_k^* +\I C_k\Lam_k^*(\a^*)^{-1}\big) \a^* S_k^{-1}\a  =0.
\label{R22}
\end{align}
The equality \eqref{R22} is immediate from \eqref{R3}, and so \eqref{R17} is also proved. Thus, \eqref{R12} is proved
as well.

Next, \eqref{R11} is proved by induction. Clearly \eqref{R9} yields \eqref{R11} for $k=0$. If \eqref{R11} holds for $k=r$,
using \eqref{R11} for $k=r$ and relations
\eqref{R1}, \eqref{R1'} and \eqref{R12} we write
\begin{align} \nn
\wt w(r+1, z)&=\left(I_m+ \frac{\I}{ z}
\wt C_r\right)\wt  w(r, z)
\\ \nn &
=\left(I_m+ \frac{\I}{ z}
\wt C_r\right)w_{\a}(r,-z)w(r,z)w_{\a}(0,-z)^{-1} 
\\ \nn &
=w_{ \alpha }(r+1, -z )\left(I_{m}+
\frac{\I}{z}C_r\right)w(r,z)w_{\a}(0,-z)^{-1}
\\  &
=w_{ \alpha }(r+1, -z )w(r+1,z)w_{\a}(0,-z)^{-1}.
\label{R23}
\end{align}
Thus, \eqref{R11} holds for $k=r+1$, and so \eqref{R11} is proved.
\end{proof}
Using \eqref{R12} we prove the next proposition.
\begin{Pn} \label{PnCk} Assume that the matrices $C_k$ satisfy the second equality in \eqref{R1}
and the triple $\{\a, S_0, \Lam_0\}$ satisfies the relations \eqref{R2}, \eqref{R6}, \eqref{R7}, and $S_0=S_0^*$.
Then, the transformed matrices $\wt C_k$ given by \eqref{R3}--\eqref{R5} have the following property$:$
\begin{align} & 
\wt C_k^2=I_m.
\label{R23+}
\end{align}

\end{Pn}
\begin{proof}. It easily follows from \eqref{R8} and \eqref{R10} (see, e.g. \cite{SaL1} or \cite[Corollary 1.13]{SaSaR}) that
\begin{align} & 
w_{\a}(r,z)w_{\a}(r,\ov{z})^*\equiv I_m.
\label{R24}
\end{align}
Since $C_k^2=I_m$, we have 
\begin{align} & 
\left(I_{m}-
\frac{\I}{z}C_k\right)
\left(I_{m}+
\frac{\I}{z}C_k\right)
=\left(1+
\frac{1}{z^2}\right)I_m.
\label{R25}
\end{align}
In view of \eqref{R24} and \eqref{R25} we derive
\begin{align} & 
w_{ \alpha }(k+1, z )\left(I_{m}-
\frac{\I}{z}C_k\right)
\left(I_{m}+
\frac{\I}{z}C_k\right)w_{ \alpha }(k+1, \ov{z} )^*
=\left(1+
\frac{1}{z^2}\right)I_m.
\label{R26}
\end{align}
On the other hand \eqref{R24} yields
\begin{align} \nn  
\left(I_{m}- \frac{\I}{z}\wt C_k\right)w_{ \alpha }(k,
z)w_{ \alpha }(k,
\ov{z})^*\left(I_{m}+ \frac{\I}{z}\wt C_k\right)&=
\left(I_{m}- \frac{\I}{z}\wt C_k\right)\left(I_{m}+ \frac{\I}{z}\wt C_k\right)
\\ & \label{R27}
=I_m+\frac{1}{z^2}\wt C_k^2.
\end{align}
According to \eqref{R12}, the left hand sides of \eqref{R26} and \eqref{R27} are equal, and so we derive
$I_m+\frac{1}{z^2}\wt C_k^2=\left(1+
\frac{1}{z^2}\right)I_m$, that is, \eqref{R23+} holds.
\end{proof}
Now, we introduce the notion of an admissible triple $\{\a, S_0, \Lam_0\}$  and show afterwards that the admissible triples
determine $S_k>0$ ($k\in \BN_0$). The definition of an admissible triple differs somewhat from the corresponding definition in \cite{FKKS},
and the proof that $S_k>0$ uses an idea from \cite{FKRS17}.

\begin{Dn}\label{DnAdm} The triple $\{\a, S_0, \Lam_0\}$ is called admissible if $\, 0, \I \, \not\in \, \s(\a)$, $S_0>0$ and  the matrix identity
\eqref{R2} is valid.
\end{Dn}

\begin{Tm} \label{PnCk+} Let an initial Dirac system \eqref{R1} and an admissible
triple $\{\a, S_0, \Lam_0\}$ be given. Then, the conditions of Theorem \ref{Fund} are satisfied.
Moreover, we have
\begin{align} & 
S_k>0 \quad (k \in \BN_0),
\label{R28}
\end{align}
and the transformed system \eqref{R1'} is skew-selfadjoint Dirac, that is, \eqref{R1+} is valid.
\end{Tm}
\begin{proof}.  In order to prove  \eqref{R28} consider the difference
\begin{align}\nn 
S_{k+1}-\big(I_n-\I \a^{-1}\big)S_k \big(I_n+\I (\a^{-1})^*\big)= &S_{k+1}-S_k-\a^{-1}S_k(\a^{-1})^*
\\ &
+\I \big(\a^{-1}S_k-S_k(\a^{-1})^*\big).
\label{R29}
\end{align}
Using \eqref{R4}, \eqref{R8} and the second equality in \eqref{R1}, we rewrite \eqref{R29} and derive a useful inequality:
\begin{align}\nn 
S_{k+1}-\big(I_n-\I \a^{-1}\big)S_k \big(I_n+\I (\a^{-1})^*\big)= & \a^{-1}\Lam_kC_k \Lam_k^*(\a^{-1})^*
\\ &
+ \a^{-1}\Lam_k \Lam_k^*(\a^{-1})^* \geq 0.
\label{R30}
\end{align}
Since $\, 0, \I \, \not\in \, \s(\a)$, the sequence $\big(I_n-\I \a^{-1}\big)^{-k}S_k \big(I_n+\I (\a^{-1})^*\big)^{-k}\,\,$ ($k\in \BN_0$) is well-defined.
In view of \eqref{R30}, this sequence is nondecreasing. Hence, taking into account $S_0>0$ we have
$\big(I_n-\I \a^{-1}\big)^{-k}S_k \big(I_n+\I (\a^{-1})^*\big)^{-k}>0$, and so \eqref{R28} holds.

Similar to \cite[Lemma A.1]{FKKS} one can show that
\begin{align} & 
\s(\a) \subset \ov{\BC_+}.
\label{R31}
\end{align}
That is, one rewrites \eqref{R2} in the form
$$
 { \left(S_0^{-1/2} \a  S_0^{1/2}\right)- \left(S_0^{-1/2} \a  S_0^{-1/2}\right)^*=\I S_0^{-1/2}\Lam_0 \Lam_0^*S_0^{-1/2},    }
$$
and from $S_0^{-1/2}\Lam_0 \Lam_0^*S_0^{-1/2}\geq 0$, the relation $\s(\a)=\s(S_0^{-1/2} \a  S_0^{1/2})~\subset~\overline{ \BC_+}$ follows.
Clearly, \eqref{R31} yields $-\I \not\in \s(\a)$. Therefore, we may set $z=-\I$ in \eqref{R12} and (taking into account the second equality in \eqref{R1} and formula \eqref{R24}) 
we obtain
\begin{align} \nn
I_m+\wt C_k&=w_{ \alpha}(k+1, -\I )(I_m+C_k)  w_{ \alpha}(k, \I )^*
\\ &
=2w_{ \alpha}(k+1, -\I )U_k^* \begin{bmatrix}I_{m_1} \\ 0 \end{bmatrix} 
\begin{bmatrix}I_{m_1} & 0 \end{bmatrix} U_k w_{ \alpha}(k, \I )^*.
\label{R32}
\end{align}
In the same way, setting in \eqref{R12} $z=\I$  we obtain
\begin{align} &
I_m-\wt C_k=2w_{ \alpha}(k+1, \I )U_k^* \begin{bmatrix}0 \\ I_{m_2}  \end{bmatrix} 
\begin{bmatrix}0 & I_{m_2} \end{bmatrix} U_k w_{ \alpha}(k, - \I )^*.
\label{R33}
\end{align}
According to \eqref{R32} and \eqref{R33} the dimension of the subspace of $\wt C_k$ corresponding
to the eigenvalue $\la =-1$ is more or equal to $m_2$ and the dimension of the subspace of $\wt C_k$ corresponding
to the eigenvalue $\la =1$ is more or equal to $m_1$.  Thus, the representation \eqref{R1+} is immediate.
\end{proof} 
 \begin{Rk} \label{RkUk}
It follows from \eqref{R1+} that $I_m+\wt C_k \geq 0$ and that $I_m+\wt C_k $ has rank $m_1$. Hence, \eqref{R32} yields
\begin{align} 
 \begin{bmatrix}I_{m_1} & 0 \end{bmatrix} U_k w_{ \alpha}(k+1, -\I )^*=\breve q_k \begin{bmatrix}I_{m_1} & 0 \end{bmatrix} U_k w_{ \alpha}(k, \I )^*, 
\quad \breve q_k>0
\label{R34}
\end{align}
 for some matrix $\breve q_k$. In the same way, formulas \eqref{R1+} and \eqref{R33} imply that 
 \begin{align} 
 \begin{bmatrix}0 & I_{m_2}  \end{bmatrix} U_k w_{ \alpha}(k+1, \I )^*=\hat q_k \begin{bmatrix}0 & I_{m_2}  \end{bmatrix} U_k w_{ \alpha}(k, - \I )^*, 
\quad \hat q_k>0.
\label{R35}
\end{align}
Now, setting
 \begin{align} 
 \wt U_k=W_k:=\diag\{\breve q_k^{1/2}, \hat q_k^{1/2}\} \left[\begin{array}{l} 
 \begin{bmatrix}I_{m_1} & 0 \end{bmatrix} U_k w_{ \alpha}(k, \I )^*
 \\ 
 \\
 \begin{bmatrix}0 & I_{m_2}  \end{bmatrix} U_k w_{ \alpha}(k, - \I )^*
  \end{array} \right]
 \label{R36}
\end{align}
we provide expressions for some suitable unitary matrices $\wt U_k$ in the representations \eqref{R1+} of the matrices $\wt C_k$.

Indeed, according to the definition of $W_k$ in \eqref{R36} and to the relations \eqref{R32}-- \eqref{R35} we have
\begin{align} & 
\wt C_k=\frac{1}{2}\left((I_m+\wt C_k)-(I_m-\wt C_k)\right)=W_k^*j W_k,
\label{R37}
\end{align}
and it remains to show that $W_k$ is unitary.  In view of \eqref{R24} and \eqref{R36}, it is easy to see that $W_kW_k^*$ is a block
diagonal matrix$\, :$
\begin{align} & 
W_k W_k^*=\diag\{\breve \rho_k, \hat \rho_k\}>0.
\label{R38}
\end{align}
Hence, for $R_k>0$ from the polar decomposition $W_k=R_k V_k$ $(V_kV_k^*=I_m)$ we have $R_k^2=\diag\{\breve \rho_k, \hat \rho_k\}$
$($and, in particular, $R_k$ is block diagonal$)$. Therefore, \eqref{R37} may be rewritten in the form
$$\wt C_k=V_k^*\, \diag\{\breve \rho_k, - \hat \rho_k\} \, V_k \qquad (\breve \rho_k>0, \,\, \hat \rho_k>0).$$
Comparing \eqref{R1+} with the formula above, we see that all the eigenvalues of $\breve \rho_k$ and $\hat \rho_k$ equal $1$,
that is, $R_k^2=R_k=I_m$, and so $W_k=V_k$. In other words, $W_k$ is unitary.
 \end{Rk}

\section{Explicit solutions of the corresponding \\ non-stationary systems}\label{Expl}
\setcounter{equation}{0}
Recall the equalities \eqref{R17}:
\begin{align} \label{R17'} &
\Lam_{k+1}^{*} S_{k+1}^{-1} (\a^2+I_{n})=
\Lam_k^{*} S_{k}^{-1} \a^2-\I \wt C_k \Lam_k^{*}
S_{k}^{-1} \a, \quad k\in \BN_0,
\end{align}
which are basic for the construction of explicit solutions of non-stationary systems.
Introduce the semi-infinite shift block matrix $\cls$ and diagonal block matrix $\wt C$:
\begin{align} \label{R39} &
\cls:=\begin{bmatrix}
0 & I_m & 0 & 0 & \ldots \\
0 & 0 & I_m & 0 & \ldots \\
0 & 0 & 0 & I_m & \ldots \\
\ldots & \ldots & \ldots & \ldots &\ldots
\end{bmatrix},
\qquad \wt C:=\diag\{\wt C_0, \, \wt C_1, \, \wt C_2, \, \ldots \}.
\end{align}
 The semi-infinite block column $\Psi(t)$ is given by the formula
 \begin{align} \label{R40} &
\Psi(t)=Y \E^{\I t \a}, \quad Y=\{Y_k\}_{k=0}^{\infty}, \quad Y_k=\Lam_k^{*} S_{k}^{-1}.
\end{align}
It is easy to see that equalities \eqref{R17'} and Theorem \ref{PnCk+} yield the
following result.
\begin{Tm} \label{TmNonst} Let the initial Dirac system \eqref{R1} and an admissible
triple $\{\a, S_0, \Lam_0\}$ be given. Then, the matrices $\Lam_k$, $S_k$ and $\wt C_k$
are well-defined via \eqref{R3}-- \eqref{R5} for $k \in \BN_0$. Moreover, the block vector
function $\Psi(t)$ constructed in \eqref{R40} satisfies the non-stationary system
\begin{align} \label{R41} &
(I-\cls)\Psi^{\prime \prime}+\wt C \Psi^{\prime}+\cls \Psi=0,
\end{align}
where $I$ is the identity operator and $\Psi^{\prime}=\frac{d}{d t} \Psi$.
\end{Tm}

\bigskip

\noindent{\bf Acknowledgments.}
 {This research   was supported by the
Austrian Science Fund (FWF) under Grant  No. P29177.}

\begin{flushright}

A.L. Sakhnovich,\\
Fakult\"at f\"ur Mathematik, Universit\"at Wien, \\
Oskar-Morgenstern-Platz 1, A-1090 Vienna, Austria\\
e-mail: oleksandr.sakhnovych@univie.ac.at
\end{flushright}



\begin{thebibliography}{AGKS}

\bibitem{Ci}
J.L.~Cieslinski,  
Algebraic construction of the Darboux matrix revisited, {\it J. Phys. A} \textbf{42}  (2009) Paper 404003.

 \bibitem{D}
  {P.A.~Deift}, 
 {Applications of a commutation formula}, {\it Duke Math. J.}   {\bf 45} (1978) 267--310.



\bibitem{FKKS}  
B. Fritzsche, M.A. Kaashoek, B. Kirstein and A.L.~Sakhnovich,  {Skew-selfadjoint Dirac systems with rational rectangular Weyl functions: explicit solutions of direct and inverse problems and integrable wave equations},
{\it Math. Nachr.}  \textbf{289} (2016)  1792--1819.

\bibitem{FKRS17}  
B. Fritzsche,  B. Kirstein, I. Roitberg and A.L.~Sakhnovich,
 Stability of the procedure of explicit recovery of skew-selfadjoint Dirac systems from rational Weyl matrix functions, {\it Linear Algebra Appl.} {\bf 533} (2017) 428--450. 

\bibitem{FKRS}  
B. Fritzsche,  B. Kirstein, I. Roitberg and A.L.~Sakhnovich,
Continuous and discrete dynamical Schr\"odinger systems: explicit
solutions, \textit{J. Phys. A: Math. Theor.} {\bf 51} (2018) Paper   015202.

\bibitem{Ge}
 F.~Gesztesy,  {A complete spectral characterization of the
double commutation method}, {\it J. Funct. Anal.}  \textbf{117} (1993) 401--446.

\bibitem{GeT}
F.~Gesztesy and G.~Teschl, 
On the double commutation method, {\it Proc. Amer. Math. Soc.} \textbf{124} (1996) 1831--1840.



\bibitem{GKS2}
 I. Gohberg, M.A. Kaashoek and A.L. Sakhnovich,  
{Pseudocanonical systems with rational Weyl functions: explicit
formulas and applications},  {\it J. Differential Equations} \textbf{146}:2 (1998) 375--398.

\bibitem{GKS6}
  {I.~Gohberg, M.A.~Kaashoek and  A.L.~Sakhnovich},  {Scattering problems for a canonical system with a pseudo-exponential potential},
 {\it Asymptotic Analysis}   {\bf 29} (2002) 1--38.



 

\bibitem{KS}
M.A. Kaashoek and   A.L. Sakhnovich, {Discrete pseudo-canonical system and isotropic Heisenberg magnet},  {\it J. Funct. Anal.}
\textbf{228}  (2005) 207--233.	 

\bibitem{MS}
  {V.B.~Matveev and M.A.~Salle},   {\it Darboux transformations and solitons},  Springer, Berlin,  1991.

\bibitem{SaA1}
  {A.L.~Sakhnovich},  
  {Exact solutions of nonlinear equations
and the method of operator identities}, {\it Linear  Alg. Appl.}    {\bf 182} (1993) 109--126.


\bibitem{SaA17}   
A.L. Sakhnovich, {Dynamics of electrons and explicit
solutions of Dirac--Weyl systems}, \textit{J. Phys. A: Math. Theor.} {\bf 50} (2017) Paper 115201.


 \bibitem{SaSaR}
 A.L.~Sakhnovich,  L.A.~Sakhnovich, and I.Ya.~Roitberg,    \textit{Inverse problems and nonlinear evolution equations. 
 Solutions, Darboux matrices and Weyl--Titchmarsh functions}, {De Gruyter Studies in Mathematics} \textbf{47},  De Gruyter, Berlin, 2013.

\bibitem{SaL1}
 {L.A. Sakhnovich}, {On  the  factorization  of  the 
transfer matrix
function},  {\it Sov. Math. Dokl.} \textbf{17} (1976) 203--207.


\bibitem{SaL2}
{L.A. Sakhnovich,} \textit{Factorisation  problems  and 
operator
identities,}   Russian
 Math. Surv. \textbf{41} (1986), 1--64.

\bibitem{ZM}
 {V.E.~Zakharov,  A.V.~Mikhailov}, 
 {On the integrability of classical spinor models in two-dimensional space-time}, {\it Commun.
Math. Phys.}   {\bf 74} (1980) 21--40.

\end{thebibliography}
\end{document}